\documentclass[12pt]{article}
\usepackage{amsfonts, amsmath, amssymb,amsthm,bbm,comment,fancyhdr,float,makeidx,mathrsfs,verbatim}
\usepackage{latexcad}

\normalsize
\newcommand{\al}{\alpha}
\newcommand{\be}{\beta}

\newcommand{\De}{\Delta}

\newcommand{\Ga}{\Gamma}

\newcommand{\seq}{\subseteq}

\newcommand{\eset}{\emptyset}

\newcommand{\ol}{\overline}

\newcommand{\vs}{\vspace*}

\newcommand{\nin}{\noindent}

\newtheorem{mthm}{Theorem}[section]
\newtheorem{mylem}[mthm]{Lemma}
\newtheorem{myprn}[mthm]{Proposition}
\newtheorem{mycor}[mthm]{Corollary}
\newtheorem{mydef}[mthm]{Definition}
\newtheorem{myrem}[mthm]{Remark}
\newtheorem{mycon}[mthm]{Construction}
\newtheorem{myeg} [mthm]{Example}
\newtheorem{myque} [mthm]{Question}

\def \nin {\noindent}

\def \Lemma #1 {\vs{3mm}\nin {\bf Lemma #1} \it}
\def \Prop #1 {\vs{3mm}\nin {\bf Proposition #1} \it}
\def \Th #1 {\vs{3mm}\nin {\bf Theorem #1} \it}
\def \Cor #1 {\vs{3mm}\nin {\bf Corollary #1} \it}
\def \Ex #1 {\vs{3mm}\nin {\bf Example #1} \it}

\def \part #1 {\hfil\break\hglue 12pt {\rm (#1)~}}

\def \qed {~\vrule height6pt width 6pt depth 0pt}
\def\fs{\footnotesize}

\title{
\bf\LARGE   Boolean graphs are unmixed and vertex decomposable\thanks{This research is supported by the National Natural
Science Foundation of China (Grant No. 11271250). } }
\author{{A-Ming Liu\thanks{aming8809@163.com} and Tongsuo Wu\thanks{Corresponding author. tswu@sjtu.edu.cn}}\\
 {\small Department of Mathematics, Shanghai Jiaotong University}
}

\date{}
\begin{document}
\baselineskip=16pt \maketitle

\begin{center}
\begin{minipage}{12cm}

 \vs{3mm}\nin{\small\bf Abstract.} {\fs For each Boolean graph $B_n$, it is proved that both $B_n$ and its complement graph $\ol{B_n}$ are vertex decomposable. It is also proved that $B_n$ is an unmixed graph, thus it is also Cohen-Macaulay.}

\vs{3mm}\nin {\small Key Words:} {\small  Boolean graph, vertex decomposable, unmixed, generalized  blow up graph }

\vs{3mm}\nin {\small 2010 AMS Classification:} {\small   Primary: 13H10; 05E45; Secondary: 13F55; 05E40.}

\end{minipage}
\end{center}

Throughout, let $[n]=\{1,2,\ldots,n\}$ and let $2^{[n]}$ be the power set of $[n]$. Recall from \cite{LUWU} that a finite Boolean graph, denoted by $B_n$, is a graph defined on the vertex set $2^{[n]}\setminus \{[n],\,\eset\}$, with $M$ adjacent to $N$ if $M\cap N=\eset$; see also \cite{LaGrange}. In the following, we list the graphs $B_3 $ and $B_4$ in diagrams (note that the center crossed in $B_4$ is not a vertex):

 \setlength{\unitlength}{0.11cm}
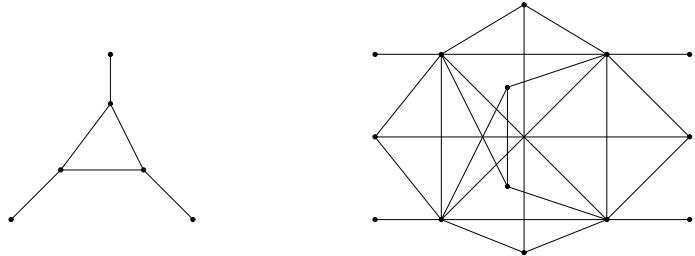
\begin{center}
\begin{picture}(94,42)
\thinlines
\drawdot{18.0}{30.0}
\drawdot{18.0}{24.0}
\drawdot{12.0}{16.0}
\drawdot{22.0}{16.0}
\drawdot{28.0}{10.0}
\drawdot{6.0}{10.0}
\drawpath{18.0}{30.0}{18.0}{24.0}
\drawpath{18.0}{24.0}{22.0}{16.0}
\drawpath{22.0}{16.0}{28.0}{10.0}
\drawpath{22.0}{16.0}{12.0}{16.0}
\drawpath{12.0}{16.0}{18.0}{24.0}
\drawpath{12.0}{16.0}{6.0}{10.0}
\drawdot{58.0}{30.0}
\drawpath{50.0}{30.0}{50.0}{30.0}
\drawpath{50.0}{30.0}{58.0}{30.0}
\drawdot{50.0}{30.0}
\drawdot{58.0}{10.0}
\drawpath{58.0}{30.0}{58.0}{10.0}
\drawpath{58.0}{10.0}{50.0}{20.0}
\drawdot{50.0}{10.0}
\drawdot{50.0}{20.0}
\drawdot{68.0}{6.0}
\drawdot{68.0}{36.0}
\drawdot{78.0}{30.0}
\drawdot{78.0}{10.0}
\drawdot{88.0}{30.0}
\drawdot{88.0}{10.0}
\drawpath{50.0}{20.0}{58.0}{30.0}
\drawpath{58.0}{30.0}{68.0}{36.0}
\drawpath{68.0}{36.0}{78.0}{30.0}
\drawpath{78.0}{30.0}{88.0}{30.0}
\drawpath{58.0}{30.0}{78.0}{30.0}
\drawpath{78.0}{30.0}{78.0}{10.0}
\drawpath{78.0}{10.0}{88.0}{10.0}
\drawpath{50.0}{10.0}{58.0}{10.0}
\drawpath{58.0}{10.0}{68.0}{6.0}
\drawpath{68.0}{6.0}{68.0}{6.0}
\drawpath{68.0}{6.0}{78.0}{10.0}
\drawpath{78.0}{10.0}{58.0}{10.0}
\drawpath{58.0}{10.0}{78.0}{30.0}
\drawpath{58.0}{30.0}{78.0}{10.0}
\drawdot{88.0}{20.0}
\drawdot{66.0}{26.0}
\drawdot{66.0}{14.0}
\drawpath{78.0}{30.0}{88.0}{20.0}
\drawpath{88.0}{20.0}{78.0}{10.0}
\drawpath{78.0}{30.0}{66.0}{26.0}
\drawpath{66.0}{26.0}{58.0}{10.0}
\drawpath{58.0}{30.0}{66.0}{14.0}
\drawpath{66.0}{14.0}{78.0}{10.0}
\drawpath{66.0}{26.0}{66.0}{14.0}
\drawpath{68.0}{36.0}{68.0}{6.0}
\drawpath{50.0}{20.0}{88.0}{20.0}
\end{picture}
\end{center}

Throughout, we use $421$ to denote the vertex $\{1,2,4\}$ of $V(B_4)$. We assume $$n>n-1>\cdots>2>1,$$ and use the pure lexicographic order on the vertices of $V(B_n)$, eg., $5421>5321$ in $V(B_6)$.

The original purpose of this work is a try to study the combinatorial property of the finite Boolean graph $B_n$, such as shellability or Cohen-Macaulayness. In the process, we find that both $B_n$ and its complement graph $\ol{B_n}$ have nice properties, and so is the related (pure) skeleton complexes and the related Alexander dual complex.

This paper is organized as follows. In section 1, we recall some basic concepts, facts and backgrounds from combinatorial commutative algebra. In section two, we first prove that $B_n$ is an unmixed graph, and then give a complicated algorithm to check that $B_n$ is also vertex decomposable. In section 3, we study the properties of the complement graph $\ol{B_n}$. In section 4, we have a preliminary study on the unmixed property of a blow up of a Boolean graph.

\section{Preliminaries}

In this part, we recall some definitions and results in combinatorial commutative algebra. For more details without mention, one can refer to the recent monographs, e.g., \cite{Villareal,Herzog and Hibi}.

Recall that a simplicial complex $\De$ is a subset of the power set $2^{[n]}$ of $[\,n\,]$, such that $\De$ is hereditary and, all singletons $x$ ($1\le x\le n$) are in $\De$. $x$ is  called a vertex  of the complex $\De$. Recall that
$$\De\setminus x=\{F\in \Delta\mid x\not\in F\},\,\,\,lk_{\De}(x)=\{F\in \Delta\mid x\not\in F,\, F\cup \{x\}\in \Delta\}.$$

\begin{mydef}  Let $\Delta$ be a simplicial complex over $[\,n\,]$. If one of the  following inductive condition is satisfied, then $\Delta$ is called
{\bf vertex decomposable:}

$(1)$ $\Delta$ is a simplex, or

$(2)$ There is a vertex $x$ such that the following requirements are fulfilled£º

\hspace{0.5cm} $(\alpha)$ Both $\Delta\setminus x$ and $lk_{\Delta}(x)$ are  vertex decomposable.

\hspace{0.5cm} $(\beta)$ No facet of $lk_{\Delta}(x)$ is a facet of $\Delta\setminus x$, or equivalently,
$$\Delta\setminus x=\langle\,\{F\mid x\not\in F\in \mathcal{F}(\Delta)\}\,\rangle.$$

\end{mydef}

Such a vertex $x$ satisfying conditions $(\al)$ and $(\be)$ is called a {\bf shedding vertex} of $\De$.  If $x$ only satisfies the second condition, then we call it a {\it weak shedding vertex}.

Recall the following implications for nonpure simplicial complexes:
$$shifted\Longrightarrow vertex\,\, decomposable\Longrightarrow shellable\Longleftarrow strongly \,\,shellable$$
Recall the following implications for simplicial complexes:
$$matroid \Longrightarrow vertex\,\, decomposable,\, and\, pure\Longrightarrow pure \,shellable$$
$$\Longrightarrow constructible\Longrightarrow Cohen-Macaulay\Longrightarrow pure.$$
For the definition of strongly shellable, see \cite{GuoWuShen}. By \cite{GuoWuShen}, if $\De$ is strongly shellable, then both $I_{\De^{\vee}}$ and $I(\De)$ have linear quotients, where $\De^\vee=\{[n]\setminus F\mid F\not\in \De\}$ and  is called the {\it Alexander dual complex} of $\De$. Note that in \cite{GuoWuShen}, counterexamples are given to show that {\it there is no implication between the concepts vertex decomposable and strongly shellable}.

The following result is well-known. Note that a similar result holds true for each of the following properties: shifted, strongly shellable, shellable,  Cohen-Macaulay, sequentially Cohen-Macaulay.

\begin{myprn}  Let $\De_1$ and $\De_2$ be complexes over $[n]=[1,n]$ and $[n+1,n+m]$ respectively. Then the join complex $\De_1*\De_2$ is vertex decomposable if and only if both complexes $\De_1$ and $\De_2$ are  vertex decomposable.
\end{myprn}

For a graph $G$, recall that the edge ideal $I(G)$ is identical with the Stanley-Reisner ideal $I_{\De_G}$ of the clique complex $\Delta_G$ of the complement graph $\ol{G}$. Recall that a graph $G$ is called {\it vertex decomposable} (Cohen-Macaulay, or shellable, or unmixed, respectively) if the simplicial complex $\De_G$ has the corresponding property. Thus we have

\begin{mycor} $(\cite[Lemma \,\,20]{Woodroofe})$\label{graph component1} A graph $G$ is vertex decomposable if and only if all connected components of $G$ are vertex decomposable.
\end{mycor}

For a vertex $x$ in a graph $G$, let $N_G[x]=N_G(x)\cup \{x\}$, the closed neighborhood of $x$ in $G$. The following is a translation of vertex decomposable of a simplicial complex in the language of a graph:

\begin{mydef} $(\cite[Lemma\, \,4]{Woodroofe})$ \label{graph version}
A graph $G$ is called vertex decomposable if either it has no edges, or else has some vertex
$x$ such that we have as follows:

$(1)$ Both $G \setminus N_G[x]$	 and $G\setminus x$ are vertex decomposable.

$(2)$ For every independent set $S$ in $G \setminus N_G[x]$, there exists some $y \in N_G(x)$ such that
$S\cup \{y\}$ is independent in $G\setminus x$.
\end{mydef}

The following result tells how to construct new and large vertex decomposable graphs from known ones:

\begin{myprn}$(\cite[Proposition 2.3]{Mousivand})$ Let $G_1,\ldots,G_n$ be finite graphs, and assume $|V(G_i)|\ge 2$, $V(G_i)\cap V(G_j)=\eset $ for all $i\not=j$. For a graph $G$ with $n$ vertices $x_i$, let $G(G_i\mid 1\le i\le n)$ be a graph obtained by attaching $x_i$ with a vertex in $G_i$. ($x_i$ is called a gluing vertex.)

 If the graphs $G_1,\ldots, G_n$ are vertex decomposable and  each gluing vertex $x_i$ is a  shedding vertex of $G_i$, then $G(G_i\mid 1\le i\le n)$ is vertex decomposable.
\end{myprn}

 Chordal graphs are an important class of vertex decomposable graphs. Recall that a graph is called {\it chordal}, if all cycles of four or more vertices have a chord, which is an edge that is not part of the cycle but connects two vertices of the cycle.  Adam Van Tuyl, Rafael H. Villarreal in \cite{VanVillareal} proved that all chordal graphs are (nonpure) shellable. Woodroofe in \cite{Woodroofe} proved further  that a chordal graph is  vertex decomposable and later, generalized the idea to clutters in \cite{Woodroofe2011}.

Recall the following  theorem, which contains important results in the algebraic combinatorics of a chordal
graph:

\begin{mthm} \label{Chordal graph in alg comb} Let $G$ be a graph and $\ol{G}$ the complement graph of $G$.
Then the following three conditions are equivalent:

$(1)$ $\ol{G}$ is chordal.

$(2)$ (Fr\"{o}berg\, \cite{Froberg}) The edge ideal $I(G)$ of $G$ has a linear
resolution.

$(3)$ (Lyubeznik\, \cite{Lyubeznik}) The cover ideal $I_c(G)$ is Cohen-Macaulay, where $I_c(G)$ is the edge ideal of the clutter consisting of all minimal vertex covers of $G$.
\end{mthm}

 Recall that a {\it simplicial vertex} of a graph  is a vertex $v$ such that the neighbourhood $N(v)$ is a clique. Recall the following main theorem (by Dirac) characterizing chordal graphs:

\begin{mthm} \label{Dirac Theorem}  A graph $G$ is chordal if and only if every induced
subgraph of $G$ has a simplicial vertex.
\end{mthm}

\section{Boolean graphs}

Recall that a {\it vertex cover} $C$ of a graph $G$ is a subset of the vertex set $V(G)$ such that
$$C\cap \{i,\,j\}\not=\varnothing,\,\,\forall \{i,\,j\}\in E(G).$$
A vertex cover is also called a {\it dominating set} of $G$, while the {\it dominating number} of $G$ is the least of cardinalities of all minimal vertex covers. Recall that a graph $G$ is said to be {\it unmixed}, if all minimal vertex covers of $G$ have a same cardinality. It is known that $G$ is unmixed if and only if the clique simplicial complex of $\ol{G}$ is unmixed, while a Cohen-Macaulay graph is always unmixed.
Recall also that $C$ is a minimal vertex cover if and only if $V(G)\setminus C$ is a maximal  independent vertex set of $G$.

Now we give the first main result of this paper:

\begin{mthm}\label{C-M} Let $n\ge 1$, and let $G$ be the Boolean graph $B_n$. Then $G$  is unmixed.
\end{mthm}

\nin {\bf Proof:} We give a  proof by considering independent vertex sets. Note  that a vertex subset  $V_0$  is a minimal vertex cover of $G$ if and only if $V_0^C$ is a maximal independent vertex set of $G$, where $V_0^C=V(G) \setminus V_0$ . In the following, we proceed to prove that all maximal independent set of $G$ have the same cardinality of $2^{n-1}-1$. In fact, a subset $\mathfrak{\Gamma}=\{b_1, b_2, \ldots, b_t\}$ of $2^{[n]}$ is  an independent vertex set of $G$, iff $b_i \bigcap b_j \neq \emptyset$ holds for distinct $b_i, b_j$ in $\mathfrak\Gamma.$  As the cardinality of the vertex set $V(G)$ is $2^{n}-2$,  we can distribute $V(G)$ into two parts $\{b_1, b_2, \ldots, b_{2^{n-1}-1}\}$ and $\{b_1^c, b_2^c, \ldots, b_{2^{n-1}-1}^c\}$. For any independent vertex set $\mathfrak{\Ga}$ of $G$ with $|\mathfrak{\Ga}|< 2^{n-1}-1$, we claim that more vertices can be added to $\mathfrak{\Ga}$ to obtain a larger independent vertex set, until the cardinality reaches $2^{n-1}-1$. For this, observe first that for any vertex $b$ in $V(G)$, the complement $b^c$ is also in $V(G)$, and this is a one to one corresponding. Thus, the cardinality of $\mathfrak{\Ga}$ is no larger than $2^{n-1}-1$. Second, for any independent vertex set $\mathfrak{\Ga}=\{b_1, \cdots b_t\}$ of $G$ with $|\mathfrak{M}|=t<2^{n-1}-1,$  clearly there exists
$b_{t+1}$ such that $\{b_{t+1},b_{t+1}^c\}\cap \mathfrak{\Ga}=\emptyset$. If neither $\mathfrak{\Ga}\cup\{b_{t+1}\}$ nor $\mathfrak{\Ga}\cup\{b_{t+1}^c\}$ is independent in the graph $G$, then there are $b_i \in \mathfrak{\Ga}$ and $b_j \in \mathfrak{\Ga}$, such that $b_i \bigcap b_{t+1}=\emptyset$ and $b_j \bigcap b_{t+1}^c=\emptyset$. Then $b_i\subseteq b_{t+1}^c$ and $b_j\subseteq b_{t+1}$, contradicting $b_i \bigcap b_j \neq \emptyset$. This shows that the graph $G$ is unmixed.
\qed

\vs{3mm} Note that $B_n$ is not chordal when $n\ge 4$, since $1-23-14-2-1$ is a cycle and it has no chord.  Note also that a Boolean graph $B_n$ is not matroidal for any $n\ge 3$. In fact, the clique complex of the complement graph $\ol{B_n}$ is far from being a  matroid in general, as the following example shows:

\begin{myeg} The clique complex of the complement graph $\ol{B_3}$ is
$$\De=\langle \{1,12,13\},\{2,12,23\},\{3,13,23\},\{12,13,23\}\rangle .$$
Note that the vertex set of $\De$ is $2^{[3]}\setminus \{[3],\emptyset\}$, so if we take a subset of it
 as $W=\{1,2,12,13\}$, then $\De_W=\langle \{1,12,13\},\{2,12\}\rangle$. Since the induced subcomplex $\De_W$ is not pure, by $\cite[Proposition\, 3.1]{Stanley}$, the complex $\De$ is not a matroid.
\end{myeg}

Next we want to prove that all Boolean graphs are vertex decomposable. In order to do so, recall that a vertex $u$ in a graph $G$ is said to have a {\it whisker}, if there is an end vertex adjacent to $u$ (\cite[Definition 7.3.10]{Villareal}). We observe the following:

\begin{mylem}\label{ whiskers} Any vertex in a graph $G$ with  whiskers is a weak shedding vertex.
\end{mylem}

\nin {\bf Proof:} Let $d$ be an end vertex adjacent to $u$. Clearly, $u\not\in G\setminus N_G[u]$, $d\in N_G(u)$,  and any independent set of $G\setminus N_G[u]$ can be extended to a larger independent vertex set $D\cup \{d\}$ in $G\setminus u$. Thus $u$ is a weak shedding vertex of $G$. \qed

\vs{3mm} Note that each of the vertices $1,\ldots,n$ has a whisker in $B_n$. If let $$G_1=B_n\setminus 1\setminus 2\setminus\cdots \setminus n,$$ then each $ji$ has a whisker in the graph $G_1$ for all $1\le i<j\le n$; If let $$G_2=G_1\setminus 12\setminus 13\setminus\cdots \setminus n-1n,$$ then every $kji$ ($n\ge k>j>i\ge 1$) has a whisker in the graph $G_2$; $\cdots$. Thus in order to show that $B_n$ is vertex decomposable, we will choose
$$1,\,\ldots,\,n;\, 21,\,\ldots,\,nn-1;\, 321,\,\ldots$$
as a sequence of weak shedding vertices. Note also that $G\setminus N_G[v]\seq G\setminus v$.

In the following, we present a weak shedding vertex order to prove the second main result of this paper:

\begin{mthm} \label{BooleanIsVERTEX DECOMPOSABLE} For any $n\ge 1$, let $G=B_n$  be the Boolean graph. Then $G$ is vertex decomposable, hence Cohen-Macaulay.
\end{mthm}

\nin {\bf Proof:} Let $$G_{n+1}=B_n, G_{n}=G_{n+1}\setminus n,G_{n-1}=G_{n}\setminus n-1,\ldots,G_{1}=G_{2}\setminus 1.\hspace{2cm} (1)$$
Note that
$G\setminus N_G[n]=\{A\cup\{n\}\mid A\in V(B_{n-1})\},$ and that it is a discrete graph, hence vertex decomposable by Corollary \ref{graph component1}.
Note also that
$$G_{n}\setminus N_{G_{n}}[n-1]=G\setminus N_{G}[n-1]\setminus n,$$
$$G_{n-1}\setminus N_{G_{n-1}}[n-2]=G\setminus N_G[n-2]\setminus n\setminus n-1,$$
$$ \ldots \ldots\ldots$$
$$ G_{2}\setminus N_{G_{n}}[1]=G\setminus N_{G}[1]\setminus n\setminus n-1\setminus \cdots\setminus 2,$$
thus they are all discrete graphs and hence, vertex decomposable. Note that each of $\{i\}$ is a weak shedding vertex of the graph $G_{i+1}$, thus by Definition \ref{graph version}, the graph $G$ is vertex decomposable if and only if the subgraph $G_{1}$ is vertex decomposable.

In order to see that the graph $G_{1}$
is vertex decomposable, let
 $$ G_{nn-1}=G_{1}\setminus n\,n-1,G_{nn-2}=G_{nn-1}\setminus nn-2,\ldots,G_{n1}=G_{n2}\setminus n1$$
 $$G_{n-1n-2}=G_{n1}\setminus \{n-1,n-2\},\,\,\ldots,\,\, G_{n-1\,1}=G_{n-1\,2}\setminus \{n-1,1\}$$
 $$\ldots\ldots\ldots\ldots,$$
 $$G_{32}=G_{41}\setminus\, 32, \,\,G_{31}=G_{32}\setminus\, 31,\,\,G_{21}=G_{31}\setminus\,\, 21 \hspace{2cm} (2)$$
Note that each vertex $ij$ is a weak shedding vertex of the graph in front of it. Now consider the corresponding $H\setminus N_H[ij]$. Let $H=G_{1}\setminus N_{G_{1}}[nn-1]$.
We have $H=(H_1\setminus nn-1)\cup (\cup_{i\ge 2} H_{2i})$, where
$$H_1=\{A\in V(B_n)\mid |A|\ge 2, n\in A\},\,H_{2i}= \{A\in V(B_n)\mid |A|=i, n-1\in A,\,n\not\in A\}.$$
Note that both $H_1$ and $\cup_{i\ge 2} H_{2i}$ are discrete graphs, and that each vertex of $H_{22}$ has a whisker  in the graph $H$ (surely, in $H_1$),
thus any linear order of vertices of $H_{22}$ is a weak shedding order of $H$. Then we delete $H_{22}$, and consider $H\setminus H_{22}$. Surely, each vertex of $H_{23}$  has a whisker in $H\setminus H_{22}$ (again, in $H_1$), thus we delete $H_{23}$ from $H\setminus H_{22}$, and continue the discussion, until we reach a forest.
This shows that $H$ is vertex decomposable.
In a similar way, we see that each of $G_{n-1i}\setminus N_{G_{n-1i}}[n-1i-1]$ is vertex decomposable.
Thus the graph $G_{1}$ is vertex decomposable if and only if the graph $G_{21}$ is vertex decomposable.

Now assume $n\ge 6$. In order to see that $G_{21}$ is vertex decomposable, the next step is to consider the sequential deletions:
$$G_{nn-1n-2}=:G_{21}\setminus \{n,n-1,n-2\}\,,\ldots,\,G_{321}=:G_{4\,2\,1}\setminus \{3,2,1\} \hspace{2cm} (3)$$
and the related  $H\setminus N_H[i\,j\,k]$.
In the process, we always take advantage of the vertices with whiskers. For the graph $L=G_{21}$, let $H=L\setminus N_L[nn-1n-2]$.
Then $$V(H)=H_1\cup(\cup_{i\ge 3}(H_{2i}\cup H_{3i})),$$ where
$$H_1=\{A\in V(B_n)\mid |A|\ge 3, n\in A,\,A\not=nn-1n-2\}$$
$$H_{2i}=\{A\in V(B_n)\mid |A|=i, n-1\in A,\,n\not\in A\}$$
$$H_{3i}=\{A\in V(B_n)\mid |A|=i, n-2\in A,\,A\cap\{ n,n-1\}=\emptyset.\}$$
Note that the subgraphs induced on each $H_i$ is discrete, and that each vertex of $H_{33}\cup H_{23}$ has a whisker in H, with an adjacent end vertex in $H_1$.
Thus in order to see that $H$ is vertex decomposable, we delete $H_{33}$ and $H_{23}$ from $H$, then going on to consider the vertices with whiskers. In this way, we show that the graph $G_{21}$ is vertex decomposable iff $G_{321}$ is vertex decomposable.

We continue this process for both related $H\setminus u$ and $H\setminus N_H[u]$, until it reaches a discrete graph or a forest.
In this way, due to the fact that the related $H\setminus N_H[u]$ always has enough weak shedding vertices (actually, vertices which have  whiskers), in the end we are able to  prove that $B_n$ is actually vertex decomposable.

Finally, it is known that vertex decomposable implies shellability, while pure shellability implies Cohen-Macaulayness. Thus by Theorem \ref{C-M}, the graph $B_n$ is Cohen-Macaulay.\quad\qed

\vs{3mm}We remark that very detailed check has been taken when $n=4,5,6$, showing that both $B_n$ and $\ol{B_n}$ are vertex decomposable. In the next section, we will show that the graph $\ol{B_n}$ is also vertex decomposable.

\section{The complement graph  $\ol{B_n}$}

Note that the graph $\ol{B_n}$ is not chordal when $n\ge 4$, since the cycle $21-32-43-41-21$ has no chord. Note also that the graph $\ol{B_n}$ is not matroidal for any $n\ge 3$. In fact, the clique complex $\De$ of $B_3$ is not pure.

Nevertheless, the complement graph also has some nice properties, see the following third main result of this paper:

\begin{mthm}\label{BooleanComplement IsVERTEX DECOMPOSABLE}  For any $n\ge 1$, the complement graph $\ol{G}$ of the Boolean graph $G=B_n$ is vertex decomposable.
\end{mthm}

\nin {\bf Proof:} For $n=3$, the result is clear. In the following, assume $n\ge 4$. Like in the Boolean case, we choose a sequence of  weak shedding vertices according to their vertex degree, and we choose it first if it has greater vertex degree. For the vertices of a same degree, we use pure lexicographic order with $n>n-1>\cdots>1$. Let
$$G_0=\ol{B_n},G_i=G_0\setminus \{\ol{1},\ldots,\ol{i},\,\}\,\,i=1,2,\ldots,n$$
where $\ol{1}=23\ldots n.$ Note that $$G_i\setminus N_{G_i}[\ol{i+1}]=\{i\},\,\,\forall\, i=0,1,\ldots,n-1,$$ and clearly condition $(2)$ of Definition \ref{graph version} is fulfilled, hence the graph $\ol{B_n}$ is vertex decomposable if and only if the graph $G_n$ is vertex decomposable.
Let
$$G_{01}=G_n,G_{12}=G_{01}\setminus \ol{12},\ldots, G_{1n}=G_{1n-1}\setminus \ol{1n},$$
$$G_{23}=G_{1n}\setminus \ol{23},G_{24}=G_{23}\setminus \ol{24},\ldots, G_{2n}=G_{2n-1}\setminus \ol{2n},$$
$$\ldots\ldots$$
$$G_{n-1n}=G_{n-2n}\setminus \ol{ n-1n}.$$
Now consider the corresponding sequence $H\setminus N_H[u]$. Note that
$$G_{01}\setminus N_{G_{01}}[\ol{12}]=\{1,2,12\}=\ol{B_2}\cup 12,$$
in which $12$ is adjacent to all vertices of $\ol{B_2}$, thus $G_{01}\setminus N_{G_{01}}[\ol{12}]$ is vertex decomposable. Since  all the corresponding $H\setminus N_H[u]$ have a same structure,
they are all vertex decomposable. Note that $3\in N_{G_{01}}(12)$, and $3$ is independent to all vertices of $G_{01}\setminus N_{G_{01}}[\ol{12}]$, thus $\ol{12}$
is a shedding vertex of the graph $G_{01}$. Similarly, it is easy to see that the sequence $$\ol{12},\ldots,\ol{1n},\,\,\ol{23},\ldots,\ol{2n},\,\,\ldots,\,\,\ol{n-1n}$$
is a shedding vertex order. Hence  the graph $G_n$ is vertex decomposable if and only if the graph $G_{n-1\,n}$ is vertex decomposable.

We continue the discussion by letting
$$G_{123}=G_{n-1n}\setminus \ol{123},G_{124}=G_{123}\setminus \ol{124}\ldots,G_{n-2\,n-1\,n}=G_{n-3\,n-1\,n}\setminus \ol{n-2\,n-1\,n}.$$
We also have
$$G_{n-1n}\setminus N_{G_{n-1n}}[\ol{123}]=2^{[3]}\setminus \ldots=\ol{B_3}\cup 123,$$
in which the vertex $123$ is adjacent to every vertex of the vertex decomposable graph $\ol{B_3}$. Note also that $4$ is shedding vertex. This shows that the graph $G_{n-1\,n}$ is vertex decomposable if and only if the graph $G_{n-2\,n-1\,n}$ is vertex decomposable. This also verifies that $\ol{B_4}$  is vertex decomposable.

If we continue this process beginning from $G_{1234}$ and ending at $G_{n-3\,n-2\,n-1\,n}$, then we proved the result for $n=5$.

This completes the verification. Clearly, this proof is a not bad algorithm, just like the proof to Theorem \ref{BooleanIsVERTEX DECOMPOSABLE}.
\qed

\vs{3mm} Recall that a {\it skeleton} complex $\De^{(0,s)}$ is a subcomplex of $\De$, which consists of all faces $F$ of $\De$ with $|F|\le s+1$. Recall that a {\it pure skeleton} complex $\De^{(s,s)}$ is a subcomplex of $\De$, which is generated by all faces of $\De$ of dimension $s$. Recall that all  skeletons and pure skeletons of a shellable complex are shellable.

By the proofs of Theorems \ref{BooleanIsVERTEX DECOMPOSABLE} and \ref{BooleanComplement IsVERTEX DECOMPOSABLE}, we have the following

\begin{mycor} Let $G$ be either the  Boolean graph $B_n$ or its complement graph $\ol{B_n}$, and let $\De$ be the clique complex of the graph $\ol{G}$. Then

$(1) $ Each skeleton complex $\De^{(0,s)}$ of $\De$ is vertex decomposable.

$(2)$  Each pure skeleton complex $\De^{(s,s)}$ of $\De$ is pure shellable, thus Cohen-Macaulay.
\end{mycor}

Note that each skeleton complex $\De^{(0,s)}$ of $\De$ is vertex decomposable if $\De$ is vertex decomposable, by \cite[Lemma 3.10]{Woodroofe2011}.

Recall that a {\it 2-flag complex} is a complex $\De$ such that each minimal nonface of
$\De$ has cardinality 2. Recall that a complex is a 2-flag complex if and only if $\De$
is a clique complex of a graph (\cite[Proposition 9.1.3]{Herzog and Hibi}). Note
that the Alexander dual $\De^\vee$ of a 2-flag complex is pure of dimension $|V(\De)|-2$.

\begin{myprn} Let $G$ be either the  Boolean graph $B_n$ or its complement graph $\ol{B_n}$, and let $\De$ be the clique complex of the graph $\ol{G}$. Then the Alexander dual complex $\De^\vee$ is not shellable when $n\geq 4$.
\end{myprn}

\nin {\bf Proof:} Recall that a complex
$\De$ is shellable if there is a shelling order of the facets $F_1, F_2, \ldots, F_t$ such that for all $i$ and $k$ with $1\leq i<k\leq t$, there exist $1\leq j<k$ and $x\in F_k$, such that $F_i\cap F_k\subseteq F_j\cap F_k=F_k\setminus\{x\}$. In the following, assume $n\ge 4.$

$(1)$ Let $\De$ be the clique complex of $\ol{B_n}$. Clearly,
$$\mathcal{F}(\De^\vee)=\{V \setminus \{a,b\}| a\in V(B_n), b\in V(B_n),\,a \cap b=\eset \} ,$$ where $V=2^{[n]}\setminus \{[n],\,\eset\}$.

Since $n\ge 4,$ we can choose $a,b\in V(B_n)$, say, $a=\{1,2\},b=\{2,3\}$, such that
$$u\cap v\neq \emptyset,\,\,\forall \,u\not=v^c, u,v\in \{a,b,a^c,b^c\}.$$
Let $F_i=V\setminus\{a,a^c\},\,F_k=V\setminus\{b,b^c\}.$

If assume that $\De^\vee$ is shellable, we can assume $F_i<F_k$ in the shelling of facets. Then by definition, there exist $1\leq j<k$ and $x\in F_k$, such that $F_i\cap F_k\subseteq F_j\cap F_k=F_k\setminus\{x\}$. If let  $F_j=V\setminus \{c,d\}$, then we have the following two facts:

$(i)$ $ F_i\cap F_k\subseteq F_j\cap F_k=F_k\setminus\{x\},$ i.e.,
$$V\setminus \{a,a^c,b,b^c\}\seq V\setminus \{c,d,b,b^c\}= V\setminus \{x,b,b^c\}.$$ It follows that $\{x,b,b^c\}=\{c,d,b,b^c\}$ and $\{c,d,b,b^c\}\subseteq \{a,a^c,b,b^c\}$

$(ii)$ $x\notin F_i$ and $x\notin F_j$, in which  $F_i=V\setminus\{a,a^c\}$ and $F_j=V\setminus \{c,d\}$.

By $(ii)$, $x\notin F_i=V\setminus \{a,a^c\}$, thus $x\in\{a,a^c\}$. Assume $x=a$, and assume further $c=a$ by fact $(i)$. Then $d\in \{b,b^c\}$ since $\{x,b,b^c\}=\{c,d,b,b^c\}$.
But then $c\cap d\not=\eset$ by the choice of $a$ and $b$,  contradicting to the assumption that $F_j$ is a facet of $\De^\vee$. The contradiction shows that $\De^\vee$ is not shellable, thus the edge ideal $I(B_n)$ does not have linear quotients.

$(2)$ As for the clique complex $\De$ of $B_n$, clearly
$$\mathcal{F}(\De^\vee)=\{V \setminus \{a,b\}|  a\in V(B_n), b\in V(B_n),\,a \cap b\neq \eset \} .$$
When $n\geq 4$, we can take $a,b,c,d\in V(B)$ with
$$a\cap b= \emptyset, \,a\cap c= \emptyset=a\cap d,\, b\cap c= \emptyset=b\cap d,$$
and  consider
$$F_i=V\setminus\{a,b\},\, F_k=V\setminus\{c,d\}.$$
If $\De^\vee$ is shellable, we can assume $F_i<F_k$ in the shelling of facets. Then a similar discussion leads to a contradiction. The details will be  omitted.
\qed

\vs{3mm} We end this section by posing the following unsettled questions:

\begin{myque} Let $G$ be either the  Boolean graph $B_n$ or its complement graph $\ol{B_n}$, and let $\De$ be the clique complex of the graph $\ol{G}$.

 $(1)$ Are the pure skeleton complexes $\De^{(s,s)}$ of $\De$ vertex decomposable?

 $(2)$ Is $\De$  strongly shellable?
 \end{myque}


\section{Blow up of Boolean graphs and unmixed property}

Recall that to get a finite  {\it blow-up graph} $G_T$ of a finite graph $G$ is to replace every vertex $v$ of $G$ by a finite set $T_v$ to get a possibly new and larger graph $G_T$, where $v\in T_v.$  The induced subgraph of $G_T$ on $T_v$ is a discrete graph, while for distinct vertices  $u,v$ of $G$, $u$ is adjacent to $v$ in $G$ if and only if each vertex of $T_u$ is adjacent to all vertices of $T_v$ in $G_T$, see \cite{KSS,NIK} for details.

If we further let $T_v$ be a complete graph, then $G_T$ becomes an {\it expanding} graph $G_E$ of $G$ (\cite{Moradi}). For a graph $G$, let $\ol{G}$ be the complement graph of $G$ in a complete graph with vertex set $V(G)$. Then the following observation holds true:

{\it A graph $H$ is a blow up of a graph $G$ if and only if $\ol{H}$ is an expanding graph of the graph $\ol{G}$.}

Note that in a non-discrete Cohen-Macaulay bipartite graph, there exists an end vertex.
Clearly, graph blow up does not keep anyone of the following properties of the original graph: chordal, vertex decomposable, Cohen-Macaulay.  On the other hand, expanding a graph keeps a lot of properties unchanged, e.g,  chordal, vertex decomposable, shellable, see \cite{Moradi} for some further discussion. Actually, for a graph, the result for chordal follows directly from Theorem \ref{Dirac Theorem}, while that for vertex decomposable follows from Definition \ref{graph version}.

In general, a  blow up of a Boolean graph is not unmixed. For example, the complete bipartite graph $K_{m,n}$ is a blow up of the Boolean graph  $B_2$ and, it is unmixed if and only if $m=n$.

\begin{myeg} Let $G_T$ be a finite blow up of the graph $B_n$. For any vertex $u\in B_n$, let $x_u=|T_u|$. Then

$(1)$ For $n=2$,  $G_T$ is unmixed if and only if $G_T=K_{m,m}$ for some $m\ge 1$.

$(2)$ For $n=3$, $G_T$ is unmixed if and only if
$x_i=x_{j\,k},\,\,\forall\,\, \{i,j,k\}=[3].$

$(3)$ For $n=4$, $G_T$ is unmixed if and only if the following seven equalities hold true:
$$x_i=x_{j\,k\,l},\,x_{i\,j}=x_{k\,l}, \,\,\forall\,\, \{i,j,k,l\}=[4].$$

$(4)$ The Boolean graphs $B_2$, $B_3$  and  $B_4$ are unmixed.
\end{myeg}

\nin {\bf Proof:}  First, note the following observations: If a graph $G$ contains a clique $K$ of $r$ vertices, then any minimal vertex cover of $G$ contains at least $r-1$ vertices of $K$; also, $G_T$ has a minimal vertex cover
which contains $\cup_{i=1}^n T_i$.

$(i)$ For $n=2$, the result is clear.

$(ii)$ For $n=3$, consider the following four
minimal vertex covers of $G_T$: $$T_{{1}}\cup T_{{2}}\cup T_{{3}},\,\,T_{{i}}\cup T_{{j}}\cup T_{\{i,j\}}(1\le i<j\leq 3).$$
Clearly, $G_S$ is unmixed if and only if the vector $(x_1,x_2, x_3,x_{11},x_{22},x_{33})$ is the positive solution in $\mathbb Z^6$ of the following system of equations:
\begin{equation}
\begin{cases}
x_1+x_2+x_{1\,2}=x_1+x_3+x_{1\,3}\\
x_1+x_2+x_{1\,2}=x_2+x_3+x_{2\,3}\\
x_1+x_2+x_{1\,2}=x_1+x_2+x_{3}\\
\end{cases}.
\end{equation}
Then the result follows.  In particular, it shows that the Boolean graph $B_3$ is unmixed.

$(iii)$ For $n=4$, note that $(\cup_{i=1}^4 T_i)\cup(\cup_{i=1}^3T_{u_i})$ is a minimal vertex cover of $G_T$, where $u_1,u_2,u_3$ are taken from distinct $\{ij,kl\}$ with  $\{i,j,k,l\}=[4]$ respectively. There are totally eight such minimal  vertex covers of $G_T$. Also, there are four others, and one representative of them is
$$(\cup_{i=2}^4T_i)\cup T_{234}\cup T_{23}\cup T_{24}\cup T_{34}.$$
Like the $n=3$ case, it follows from the system of linear equations that $x_i=x_{j\,k\,l}$ holds for all $\{i,j,k,l\}=[4]$. Then
it follows easily $x_{i\,j}=x_{k\,l}$.

The converse holds clearly.

In particular, the Boolean graph $B_i$ ($1\ge i\le 4$) is unmixed. \qed

\vs{3mm}This shows another way for illustrating Theorem 2.1.  When $n$ is large, things will become complicated. But a similar careful  discussion shows that the unmixedness of the blow up $G_T$ of the Boolean graph $B_n$ ($n=5,6,7$,  respectively) amounts to the solving of a system of linear equations with indeterminate labeled properly according to their position in the layers.

The above example shows that {\it graph blow up} is a good concept for discussing unmixed property of graphs. We can even generalize it a little  to obtain a {\it finite generalized  blow up} $G_S$ of a finite graph $G$ explained in what follows. For every  vertex $v$ of $G$, let $S_v$ be a disjoint union of $S_{1v}$ with $S_{2v}$, in which $v\in S_{1v}$. Replace $v$ by $S_v$ to get a possibly new and larger graph $G_S$: For any $u\in V(G)$, the induced subgraph of $G_S$ on each $S_u$ is a discrete graph, while for distinct vertices  $u,v$ of $G$, $u$ is adjacent to $v$ in $G$ iff each vertex of $S_{1u}$ is adjacent to all vertices of $S_v$ and each vertex of $S_{1v}$ is adjacent to all vertices of $S_u$. Note that whenever none of $ S_{2u}, S_{2v}$ is empty, no vertices in $ S_{2u} $ is adjacent to a vertex in $ S_{2v}$. By the definition, each blow up is a generalized blow up, of a graph; but the converse is clearly not true.

Generalized blow up occur naturally when we consider deleting a vertex from the graph $B_n$, as the following example shows.

\begin{myeg} $B_n\setminus n\setminus 12\ldots n-1$ is a generalized blow up of $B_{n-1}$.
\end{myeg}

\nin {\bf Proof:}   Clearly, the vertex $12\ldots n-1$ is isolated in the graph $B_n\setminus n$.

 Let $G=B_n\setminus n\setminus 12\ldots n-1.$ Then the vertex set of $V(G)$ splits with two parts, $\{ A, A\cup\{n\}\},\forall A\in V(B_{n-1})$. Thus if we add $A\cup\{n\}$ to
 the vertex $A$ as the second part, then clearly, $G$ is a generalized blow up of $B_{n-1}$, where for each vertex $v$ of $B_{n-1}$, we have $|S_{1v}|=|S_{2v}|=1.$
\qed

\vs{3mm}We end the paper with an easy discussion on the unmixedness of a generalized blow up of the graph $G=B_n.$

\begin{myeg} Let $G_S$ be a generalized blow up of the graph $G=B_2$. Then $G_S$ is unmixed if an only if either $G_S=K_{m,m}$ or $|S_{11}|=|S_{12}|=|S_{21}|=|S_{22}|.$
\end{myeg}

\nin {\bf Proof:} Assume that $G_S$ is a generalized blow up of the graph $G=B_2$, but not a blow up of $B_2$. Note that $S_{1\{1\}}\cup S_{1\{2\}}$, $S_{2\{1\}}\cup S_{2\{2\}}$  and $S_{1\{1\}}\cup S_{2\{1\}}$ are minimal vertex covers of the graph $G_S$. Thus if $G_S$ is unmixed, then we have
 $$|S_{1\{1\}}|+|S_{1\{2\}}|=|S_{2\{1\}}|+|S_{2\{2\}}|=|S_{1\{1\}}|+|S_{2\{1\}}|$$
hence all $|S_{ij}|$ are identical.

The converse
holds clearly.

\end{document}